\documentclass[twoside]{article} 
\usepackage{graphicx}
\date{} 
\title{The asymptotic expansion of the Bateman and Havelock functions of large order and argument}
\author{\sc R. B.\ Paris \\
{\em Division of Computing and Mathematics,} \\
{\em Abertay University, Dundee DD1 1HG, UK}}
\begin{document}
\def\f#1#2{\mbox{${\textstyle \frac{#1}{#2}}$}}
\def\dfrac#1#2{\displaystyle{\frac{#1}{#2}}}
\def\boldal{\mbox{\boldmath $\alpha$}}
\newcommand{\bee}{\begin{equation}}
\newcommand{\ee}{\end{equation}}
\newcommand{\sa}{\sigma}
\newcommand{\ka}{\kappa}
\newcommand{\al}{\alpha}
\newcommand{\la}{\lambda}
\newcommand{\ga}{\gamma}
\newcommand{\eps}{\epsilon}
\newcommand{\om}{\omega}
\newcommand{\fr}{\frac{1}{2}}
\newcommand{\fs}{\f{1}{2}}
\newcommand{\g}{\Gamma}
\newcommand{\br}{\biggr}
\newcommand{\bl}{\biggl}
\newcommand{\ra}{\rightarrow}
\newcommand{\gtwid}{\raisebox{-.8ex}{\mbox{$\stackrel{\textstyle >}{\sim}$}}}
\newcommand{\ltwid}{\raisebox{-.8ex}{\mbox{$\stackrel{\textstyle <}{\sim}$}}}
\renewcommand{\topfraction}{0.9}
\renewcommand{\bottomfraction}{0.9}
\renewcommand{\textfraction}{0.05}
\newcommand{\mcol}{\multicolumn}
\date{}
\maketitle
\pagestyle{myheadings}
\markboth{\hfill \sc R. B.\ Paris  \hfill}
{\hfill \sc Asymptotics of the Bateman and Havelock functions\hfill}
\begin{abstract}
Asymptotic expansions for the Bateman and Havelock functions defined respectively by the integrals
\[\frac{2}{\pi}\int_0^{\pi/2} \!\!\!\begin{array}{c} \cos\\\sin\end{array}\!(x\tan u-\nu u)\,du\]
are obtained for large real $x$ and large order $\nu>0$ when $\nu=O(|x|)$. The expansions are obtained by application of the method of steepest descents combined with an inversion process to determine the coefficients.
Numerical results are presented to illustrate the accuracy of the different expansions obtained.
\vspace{0.3cm}

\noindent {\bf Mathematics subject classification (2020):} 30E20, 33E20, 34E05, 41A60 
\vspace{0.1cm}
 
\noindent {\bf Keywords:} Bateman function, Havelock function, asymptotic expansions, method of steepest descents
\end{abstract}

\vspace{0.3cm}

\noindent $\,$\hrulefill $\,$

\vspace{0.3cm}

\begin{center}
{\bf 1.\ Introduction}
\end{center}
\setcounter{section}{1}
\setcounter{equation}{0}
\renewcommand{\theequation}{\arabic{section}.\arabic{equation}}
The Bateman function, defined by the integral
\bee\label{e11}
k_\nu(x)=\frac{2}{\pi}\int_0^{\pi/2} \cos (x\tan u-\nu u)\,du,
\ee
was introduced by Bateman in 1931 \cite{HB} as a solution of the ordinary differential equation
\[xu''(x)=(x-\nu)u(x)\]
for non-negative integer values of $\nu$
arising in von K\'arm\'an's theory of turbulent flows. An analogous integral, appearing in the problem of surface waves, was introduced by Havelock \cite{THH} as
\bee\label{e12}
h_\nu(x)=\frac{2}{\pi}\int_0^{\pi/2} \sin (x\tan u-\nu u)\,du.
\ee

An excellent and detailed survey of the properties of these two functions, including certain generalisations and integrals, has recently been given by Apelblat, Consiglio and Mainardi \cite{ACM}. Expressions for $k_n(x)$ when $n=0, 1, 2, \ldots$ are listed in \cite[\S\S2, 3]{ACM} and involve $e^{-x}$ multiplied by a polynomial in $x$ of degree $\fs n$ for even values of $n$ and the modified Bessel function of the second kind for odd values of $n$. Expressions for $h_n(x)$ are more complicated \cite[\S 3]{ACM} and involve the logarithmic integral li$(e^x)$, where li$(z)=\int_0^{z}dt/\log\,t$.

The Bateman function can be expressed in terms of the confluent hypergeometric function $U(a,b,z)$, which is given in
\cite[p.~510]{AS} as
\bee\label{e13}
k_\nu(x)=\frac{e^{-x}}{\g(1+\fs \nu)}\,U(-\fs\nu,0,2x)= \frac{1}{\g(1+\fs\nu)}\,W_{\nu/2,1/2}(2x) \qquad (x>0),\
\ee
where $W_{\kappa,\mu}(z)$ is the Whittaker function \cite[p.~334]{DLMF}. In the case of negative values of $x$, we have
\begin{eqnarray}
k_\nu(-x)\!\!&=&\!\!\frac{2}{\pi}\int_0^{\pi/2} \cos (x\tan u+\nu u)\,du\nonumber\\
&=&\!\!e^{-x} \g(\fs\nu) \sin \fs\pi\nu\, U(\fs\nu,0,2x)\qquad (x>0),\label{e14}
\end{eqnarray}
from which it is seen that $k_{2n}(-x)=0$, $n=0, 1, 2, \ldots\  $. 

The asymptotic expansions of $k_\nu(x)$ for $x\to\pm\infty$ and fixed order $\nu$ follow from (\ref{e13}) and (\ref{e14}) by using the well-known expansion for $U(a,b,z)$ \cite[(13.7.10)]{DLMF} to find
\bee\label{e15a}
k_\nu(x)\sim\frac{(2x)^{\nu/2}e^{-x}}{\g(1+\nu/2)} \sum_{k=0}^\infty \frac{(-)^k (-\fs\nu)_k (1-\fs\nu)_k}{k! (2x)^k} \qquad(x\to+\infty)\ee
\bee\label{15b}k_\nu(-x)\sim(2x)^{-\nu/2}e^{-x} \g(\fs\nu) \sin \fs\pi\nu\,\sum_{k=0}^\infty \frac{(-)^k (\fs\nu)_k (1+\fs\nu)_k}{k! (2x)^k}\qquad (x\to+\infty),
\ee
where $(a)_k=\g(a+k)/\g(a)$ is the Pochhammer symbol. The expansions of the Havelock function for $x\to\pm\infty$ and fixed $\nu$ are given by (for details see Appendix A)
\bee\label{e16}
h_\nu(x)\sim \frac{2}{\pi x} \sum_{k=0}^\infty \frac{c_k(\nu)}{x^k},\qquad h_\nu(-x)\sim\frac{2}{\pi x} \sum_{k=0}^\infty \frac{(-)^{k-1}c_k(\nu)}{x^k}
\ee
as $x\to+\infty$, where an exponentially small contribution of $O(e^{-x})$ has been neglected in each case. The first two coefficients are  $c_0(\nu)=1$, $c_1(\nu)=\nu$; the coefficients $c_k(\nu)$ are listed in the appendix for $k\leq 8$.
It is seen that the two functions have completely different asymptotic behaviours: the Bateman function is
exponentially decaying as $x\to\pm\infty$ whereas the Havelock function possesses an algebraic decay of O$(x^{-1})$.

In this paper we consider the asymptotic expansions of $k_\nu(x)$ and $h_\nu(x)$ of large argument and order. We shall suppose throughout that $x$ and $\nu$ are real variables and, in the case $x>0$, set $\nu=ax$, where $a>0$ is a fixed parameter. In the case of $k_\nu(x)$, it would, of course,  be possible to exploit the known asymptotics of the Whittaker function $W_{\kappa,\mu}(x)$ for large $\kappa>0$ and finite $\mu$ uniformly valid in $x\in[\delta,\infty)$ \cite[p.~412, Ex.~7.3]{Olv}; see also \cite[\S 13.21]{DLMF}. However, the resulting expansion and coefficients are complicated and we find it easier here to derive the asymptotics of $k_\nu(x)$ and $h_\nu(x)$, and hence also 
$k_\nu(\nu)$ and $h_\nu(\nu)$, directly by the method of steepest descents applied to the integrals (\ref{e11}) and (\ref{e12}).
In Sections 2 and 3, we present the details for the case $x\to+\infty$ with $\nu=O(x)$ and in Section 4 the case $x\to-\infty$ with $\nu=O(|x|)$. In Section 5 we display numerical results confirming the accuracy of the different expansions obtained. 

\vspace{0.6cm}

\begin{center}
{\bf 2.\ The asymptotic expansion of $k_\nu(x)$ when $x\to+\infty$, $\nu=ax$}
\end{center}
\setcounter{section}{2}
\setcounter{equation}{0}
\renewcommand{\theequation}{\arabic{section}.\arabic{equation}}
We first consider the Bateman function with $\nu=ax$, which from (\ref{e11}) we write as
\bee\label{e21}
k_\nu(x)=\frac{1}{\pi}\int_{-\pi/2}^{\pi/2} \cos (x(\tan u-au))\,du=\frac{1}{2\pi}\int_{-\pi/2}^{\pi/2} \{e^{x\psi(u)}+e^{-x\psi(u)}\}\,du,
\ee
where
\bee\label{e21a}
\psi(u)=i\phi(u)=i(\tan u-au).
\ee

Saddle points of $\psi(u)$ occur when $\psi'(u)=i(\mbox{sec}^2u-a)=0$ and are given by
\[u_0=\pm \arccos \frac{1}{\sqrt{a}}=\pm\arctan \sqrt{a-1}.\]
When $a>1$, the saddles are situated on the real $u$-axis in $[-\fs\pi,\fs\pi]$, and when $a<1$ they are situated on the imaginary axis at $u_0=\pm i\,\mbox{arctanh} \sqrt{1-a}$. In the case $a=1$ the saddles coalesce to form a double saddle point at the origin.

\vspace{0.3cm}

\noindent{\bf 2.1}\ \ {\bf The case $a>1$}
\vspace{0.3cm}

\noindent The paths of steepest descent through the saddles when $a>1$ are shown in Fig.~1(a). Since $\psi(iy)=ay-\tanh y$, it is seen that $\psi(iy)\to\pm\infty$ as $y\to\pm\infty$. Consequently, the integration path for the part of the integrand in (\ref{e21}) involving $\exp\{x\psi(u)\}$ can be deformed to pass over the paths labelled $ABCD$, whereas that part involving $\exp\{-x\psi(u)\}$ passes over the mirror image $EFGH$. It is clear that the contribution from the saddle $-u_0$ is equal to that from the saddle $u_0$, so that it is sufficient to consider only the contribution from the positive saddle. It is to be remarked that the paths at $B$ and $C$ (and at $F$ and $G$) do not asymptote to the imaginary axis, but to the lines $\Re (u)=\pm c/a$, where $c=a\,\arctan\,\sqrt{a-1}-\sqrt{a-1}>0$ when $a>1$; see Appendix B. The steepest descent paths are connected by the horizontal lines $BC$ (and $FG$), which yield no contribution as $\Im (u)\to\mp\infty$.
\begin{figure}[th]
	\begin{center}	{\tiny($a$)}\ \includegraphics[width=0.35\textwidth]{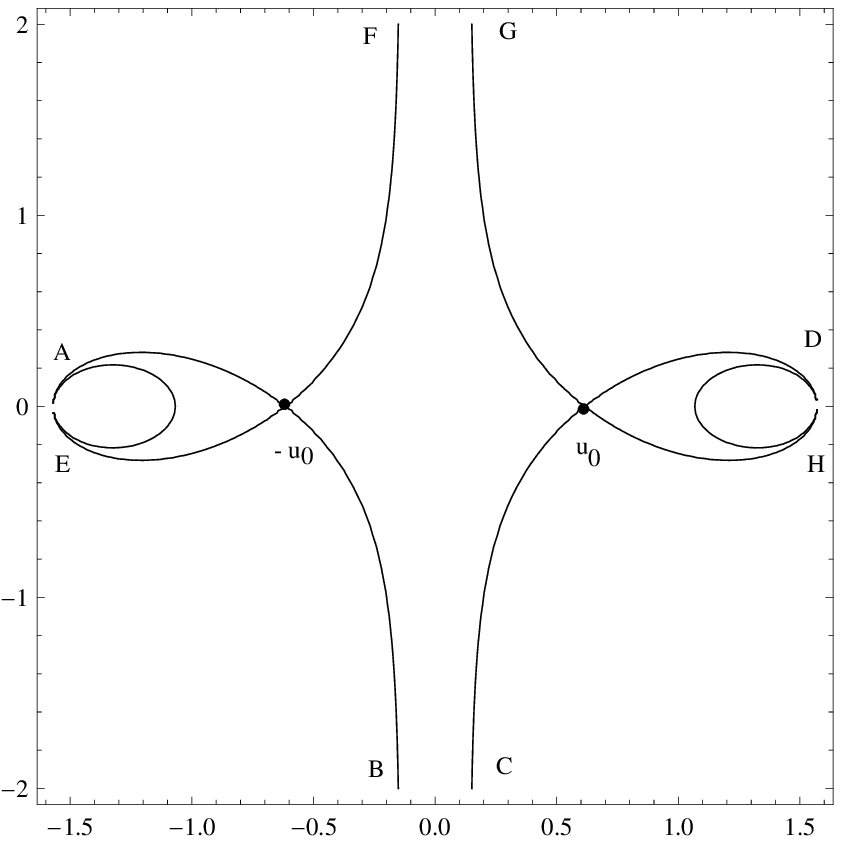}
	\qquad
	{\tiny($b$)}\ \includegraphics[width=0.35\textwidth]{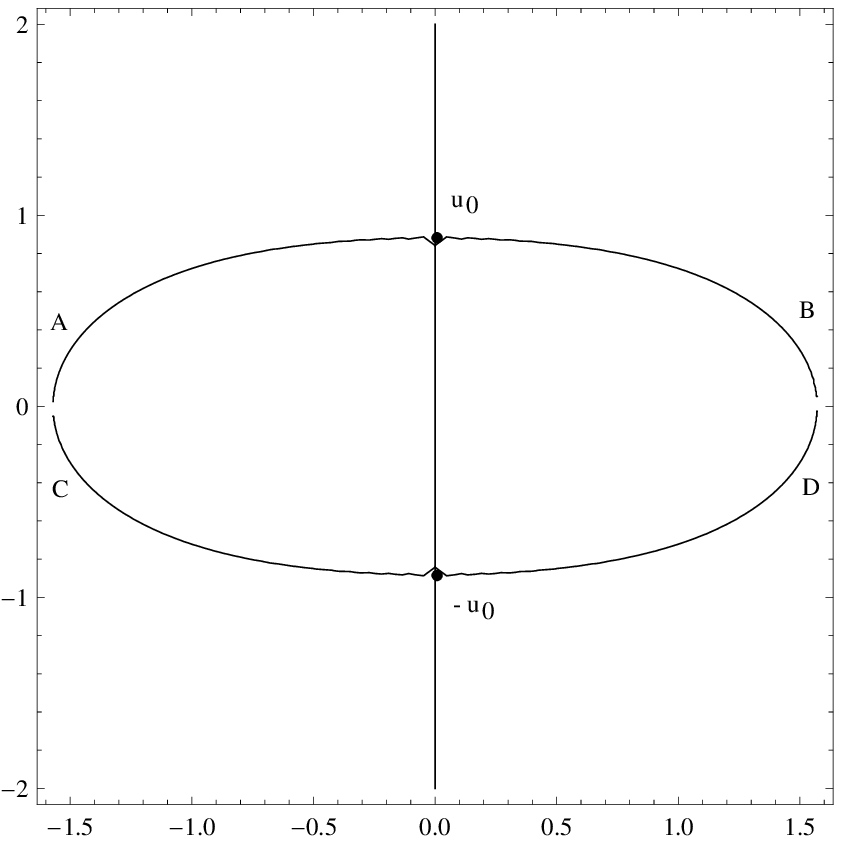} 
	\vspace{0.5cm}
	
	\qquad
	{\tiny($c$)}\ \includegraphics[width=0.35\textwidth]{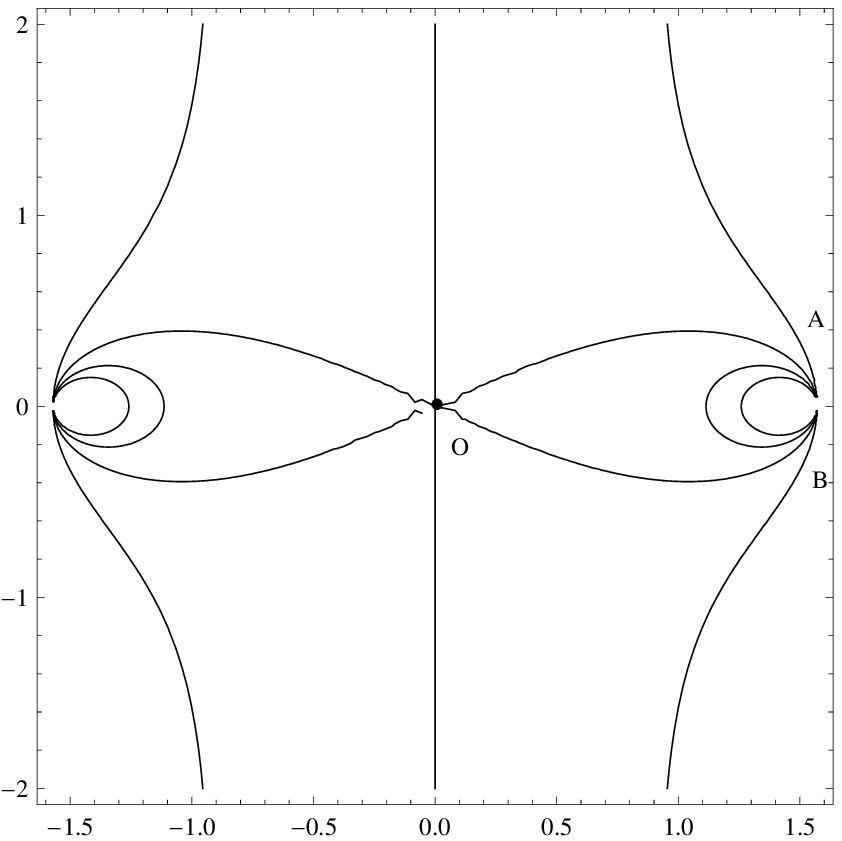}
\caption{\small{Paths of steepest descent through the saddles $\pm u_0$: (a) $a>1$, (b) $a<1$ and (c) $a=1$.}}
\end{center}
\end{figure}

At the endpoints of the steepest path $CD$ we have $\psi(u)\to-\infty$. Then
\[I_{CD}:=\frac{1}{2\pi}\int_{CD}e^{x\psi(u)}du=\frac{e^{x\psi(u_0)}}{2\pi} \int_{-\infty}^\infty e^{-xw^2} \frac{du}{dw}\,dw,\]
where 
\[\psi(u_0)=i\phi(u_0)=i(\sqrt{a-1}-a \arctan\sqrt{a-1}).\]
We have introduced the new variable $w$ by
\[-w^2=\psi(u)-\psi(u_0)=\frac{\psi''(u_0)}{2!} (u-u_0)^2+\frac{\psi'''(u_0)}{3!} (u-u_0)^3
%+\frac{\psi''''(u_0)}{4!} (u-u_0)^4
+\cdots\]
\[=ia\sqrt{a-1} (u-u_0)^2+\frac{1}{3}ia(3a-2)(u-u_0)^3+\frac{1}{3}ia\sqrt{a-1}(3a-1)(u-u_0)^4+\cdots\ .\]
Inversion of this series (essentially using Lagrange's theorem) followed by differentiation then yields
\[\frac{du}{dw}\stackrel{e}{=} \sqrt{\frac{i}{a}}(a-1)^{-1/4} \sum_{k\geq 0}i^k A_{2k}(a) w^{2k},\]
where $\stackrel{e}{=}$ signifies the inclusion of {\em only the even powers of $w$}, since odd powers will not enter into this calculation. With the help of {\it Mathematica}, using the InverseSeries command, the first few coefficients\footnote{The coefficients $A_{2}(a)$ and $A_4(a)$ can be obtained alternatively from the standard representation of the saddle-point coefficients; see, for example, \cite[pp.~13--14]{PHad}.} $A_{2k}(a)$ are given by
\[A_0(a)=1,\quad A_2(a)=\frac{(8-12a+9a^2)}{24a(a-1)^{3/2}},\quad A_4(a)=\frac{(64-192a+288a^2+360a^3-135a^4)}{3456a^2(a-1)^3},\]
\[A_6(a)=\frac{(-71168+320256a-554688a^2+518400 a^3 + 340200 a^4 - 170100 a^5 + 42525 a^6)}{6220800a^3(a-1)^{9/2}},\]
\[A_8(a)=\frac{1}{4180377600a^4(a-1)^6}(-2338816 + 14032896 a - 36790272 a^2 + 55710720 a^3 - 32876928 a^4 \]
\bee\label{e23b}+ 
 231880320 a^5 - 68584320 a^6 + 30618000 a^7 - 5740875 a^8).\ee
Then we obtain
\[I_{CD}\sim \frac{e^{i(x\phi(u_0)+\pi/4)}}{\pi \sqrt{x}}\,(a-1)^{-1/4}\sum_{k=0}^\infty i^k A_{2k}(a) \! \int_0^\infty e^{-xw^2}  w^{2k}dw\]
\bee\label{e22}
=\frac{e^{i(x\phi(u_0)+\pi/4)}}{2\sqrt{\pi ax}}\,(a-1)^{-1/4}\sum_{k=0}^\infty \frac{i^k \g(k+\fs)}{\g(\fs)}\,\frac{A_{2k}(a)}{x^k}.\hspace{0.2cm}
\ee
 
The contribution from the path $GH$ involving $\exp \{-x\psi(u)\}$ is given by the conjugate of (\ref{e22}). Taking into account that the contribution from the saddle at $u=-u_0$ is the same as that from $u=u_0$, then  after some routine algebra we finally obtain the expansion:
\newtheorem{theorem}{Theorem}
\begin{theorem}$\!\!\!.$\ \  The asymptotic expansion of $k_\nu(x)$ for $x\to+\infty$ with $\nu=ax$, $a>1$ is
\bee\label{e23}
k_\nu(x)\sim \frac{2(a-1)^{-1/4}}{\sqrt{\pi ax}}\bl\{\cos \Phi(x)\sum_{k=0}^\infty \frac{(-)^k A_{4k}(a) (\fs)_{2k}}{x^{2k}} -\frac{\sin \Phi(x)}{2x}
\sum_{k=0}^\infty \frac{(-)^k A_{4k+2}(a) (\f{3}{2})_{2k}}{x^{2k}}\br\},
\ee
where 
\bee\label{e23a}
\Phi(x):=x\phi(u_0)+\frac{\pi}{4}=x(\sqrt{a-1}-a\arctan\sqrt{a-1})+\frac{\pi}{4}
\ee
and the coefficients $A_{2k}(a)$ are defined in (\ref{e23b}). 
\end{theorem}
\vspace{0.3cm}

\noindent{\bf 2.2}\ \ {\bf The case $a<1$}
\vspace{0.3cm}

\noindent When $a<1$ the saddle points are situated on the imaginary axis at $u=\pm i\,\mbox{arctanh} \sqrt{1-a}$. The steepest descent paths through the saddles are illustrated in Fig.~1(b). The integration path for the part of the integrand in (\ref{e21}) involving $\exp\{x\psi(u)\}$ can be deformed to pass over the path labelled $AB$, whereas that part involving $\exp\{-x\psi(u)\}$ passes over the mirror image $CD$. 
 
Then, as before,
\[I_{AB}=\frac{1}{2\pi}\int_{AB}e^{x\psi(u)}du=\frac{e^{x\psi(u_0)}}{2\pi}\int_{-\infty}^\infty e^{-xw^2}\frac{du}{dw}\,dw,\]
where the new variable $w$ is given by
\[-w^2=\psi(u)-\psi(u_0)=-a\sqrt{1-a} (u-u_0)^2+\frac{1}{3}ia(3a-2)(u-u_0)^3%-\frac{1}{3}a\sqrt{1-a}(3a-1)(u-u_0)^4
+\cdots\ .\]
Upon inversion followed by differentiation, we obtain
\bee\label{e23c}
\frac{du}{dw}\stackrel{e}{=}\frac{(1-a)^{-1/4}}{\sqrt{a}} \sum_{k\geq0} {\hat A}_{2k}(a) w^{2k},
\ee
where
\[{\hat A}_0(a)=1,\quad {\hat A}_2(a)=\frac{-(8-12a+9a^2)}{24a(1-a)^{3/2}},\quad {\hat A}_4(a)=\frac{(64-192a+288a^2+360a^3-135a^4)}{3456a^2(1-a)^3},\]
\[{\hat A}_6(a)=\frac{(71168-320256a+554688a^2-518400 a^3 - 340200 a^4 + 170100 a^5 - 42525 a^6)}{6220800a^3(1-a)^{9/2}},\]
\[{\hat A}_8(a)=\frac{1}{4180377600a^4(1-a)^6}(-2338816 + 14032896 a - 36790272 a^2 + 55710720 a^3 - 32876928 a^4 \]
\bee\label{e23d}
+ 231880320 a^5 - 68584320 a^6 + 30618000 a^7 - 5740875 a^8).
\ee

Then we find
\[I_{AB}\sim \frac{e^{x\psi(u_0)}}{2\sqrt{\pi ax}}\,(1-a)^{-1/4} \sum_{k=0}^\infty\frac{\g(k+\fs)}{\g(\fs)}\, \frac{{\hat A}_{2k}(a)}{x^k},\]
where
\bee\label{e24a}
\psi(u_0)=-\Psi(a),\qquad \Psi(a):=\sqrt{1-a}-a\,\mbox{arctanh} \sqrt{1-a}.
\ee
An equal contribution results from the path $CD$ with the exponential factor $\exp \{-x\psi(u)\}$ to yield: 
\begin{theorem}$\!\!\!.$\ \ The asymptotic expansion of $k_\nu(x)$ for $x\to+\infty$ with $\nu=ax$, $a<1$ is
\bee\label{e24}
k_\nu(x) \sim \frac{(1-a)^{-1/4}}{\sqrt{\pi ax}}\,e^{-x\Psi(a)} \sum_{k=0}^\infty\frac{\g(k+\fs)}{\g(\fs)}\, \frac{{\hat A}_{2k}(a)}{x^k},
\ee
where the coefficients ${\hat A}_{2k}(a)$ are defined in (\ref{e23d}).
\end{theorem}
\vspace{0.3cm}

\noindent{\bf 2.3}\ \ {\bf The case $a=1$}
\vspace{0.3cm}

\noindent When $a=1$ the two saddles coalesce to form a double saddle at $u=0$. In this case we consider the integral (\ref{e11}) over the interval $[0,\fs\pi]$. The paths of steepest descent from the origin are shown in Fig.~1(c); the upper path $OA$ is chosen for the exponential factor $\exp \{x\psi(u)\}$ with the lower (conjugate) path $OB$ for the factor $\exp \{-x\psi(u)\}$. Then
\[I_{OA}=\frac{1}{\pi}\int_{OA}e^{x\psi(u)}du=\frac{1}{\pi}\int_0^\infty e^{-xw^3} \frac{du}{dw}\,dw,\]
where, since $\psi(0)=0$, the new variable $w$ is now given by
\[-w^3=\psi(u)=\frac{iu^3}{3}+\frac{2iu^5}{15}+\frac{17iu^7}{315}+\frac{62iu^9}{2835}+\cdots\ .\]
Inversion and differentiation then yields (where no odd powers of $w$ are present)
\[\frac{du}{dw}=\sum_{k\geq0} B_k w^{2k},\]
where, with $\mu=3^{1/3}e^{\pi i/6}$,
\[B_0=\mu,\quad B_1=-\frac{6i}{5},\quad B_2=-\frac{27}{35\mu},\quad B_3=\frac{2\mu}{25},\]
\bee\label{e2dk}
B_4=\frac{1296i}{67375}, \quad B_5= \frac{9774}{284375\mu},\quad B_6=-\frac{49711\mu}{11790625},\quad B_7=-\frac{3390336i}{1861234375},\ldots\ .
\ee

Then
\[I_{OA}\sim\frac{1}{\pi}\sum_{k\geq0} B_k \int_0^\infty e^{-xw^3} w^{2k} dw=\frac{x^{-1/3}}{3\pi}\sum_{k\geq0} \frac{B_k }{x^{2k/3}}\int_0^\infty e^{-\tau}\tau^{2k/3-2/3} d\tau\] 
\bee\label{e26}
=\frac{x^{-1/3}}{3\pi} \sum_{k=0}^\infty \frac{B_k}{x^{2k/3}}\,\g(\f{2}{3}k+\f{1}{3}).
\ee
The contribution from the lower path $OB$ yields the conjugate expansion. Hence, upon observing that the coefficients $B_k$ corresponding to $k=1, 4, 7, \ldots$ are pure imaginary, we obtain:
\begin{theorem}$\!\!\!.$\ \ The asymptotic expansion of $k_\nu(\nu)$ for $\nu\to+\infty$ is
\bee\label{e2k}
k_\nu(\nu)=I_{OA}+{\overline I}_{OA}\sim \frac{2}{3\pi \nu^{1/3}}\!\!\!\!\! 
\mathop{\sum_{k\geq0}}_{\scriptstyle k\neq 1, 4, 7, \ldots}\!\!\!\!\! \frac{\Re (B_k)}{\nu^{2k/3}}\,\g(\f{2}{3}k+\f{1}{3}),
\ee
where the coefficients $B_k$ are listed in (\ref{e2dk}).
An explicit representation of the first few terms in this expansion is
\[k_\nu(\nu)\sim\frac{3^{-1/6}}{\pi \nu^{1/3}}\bl\{\g(\f{1}{3})-\frac{9\cdot 3^{1/3}}{35\,\nu^{4/3}}\g(\f{5}{3})+\frac{2}{25\,\nu^2} \g(\f{7}{3})+\frac{3258\cdot 3^{1/3}}{284375\,\nu^{10/3}} \g(\f{11}{3})\hspace{2cm}\]
\bee\label{e25}
\hspace{7cm}-\frac{49771}{11790625\,\nu^4} \g(\f{13}{3})+\cdots\br\}.
\ee 
\end{theorem}
 
\vspace{0.6cm}

\begin{center}
{\bf 3.\ The asymptotic expansion of $h_\nu(x)$ when $x\to+\infty$, $\nu=ax$}
\end{center}
\setcounter{section}{3}
\setcounter{equation}{0}
\renewcommand{\theequation}{\arabic{section}.\arabic{equation}} 
From (\ref{e12}), the Havelock function can be written as
\[h_\nu(x)=\frac{1}{\pi i}\int_0^{\pi/2} \{e^{x\psi(u)}-e^{-x\psi(u)}\}\,du,\]
where $\psi(u)$ is defined in (\ref{e21a}). 
The integration path $[0,\fs\pi]$ will be deformed to coincide with the paths of steepest descent through the saddle point $u=u_0$ (when $a>1$) and the conjugate saddles (when $a<1$). These paths situated in $\Re (u)\geq 0$ are shown in Fig.~2 and are basically the same as those in Fig.~1, but with the addition of a connecting path emanating from the origin.
\begin{figure}[th]
	\begin{center}	{\tiny($a$)}\ \includegraphics[width=0.35\textwidth]{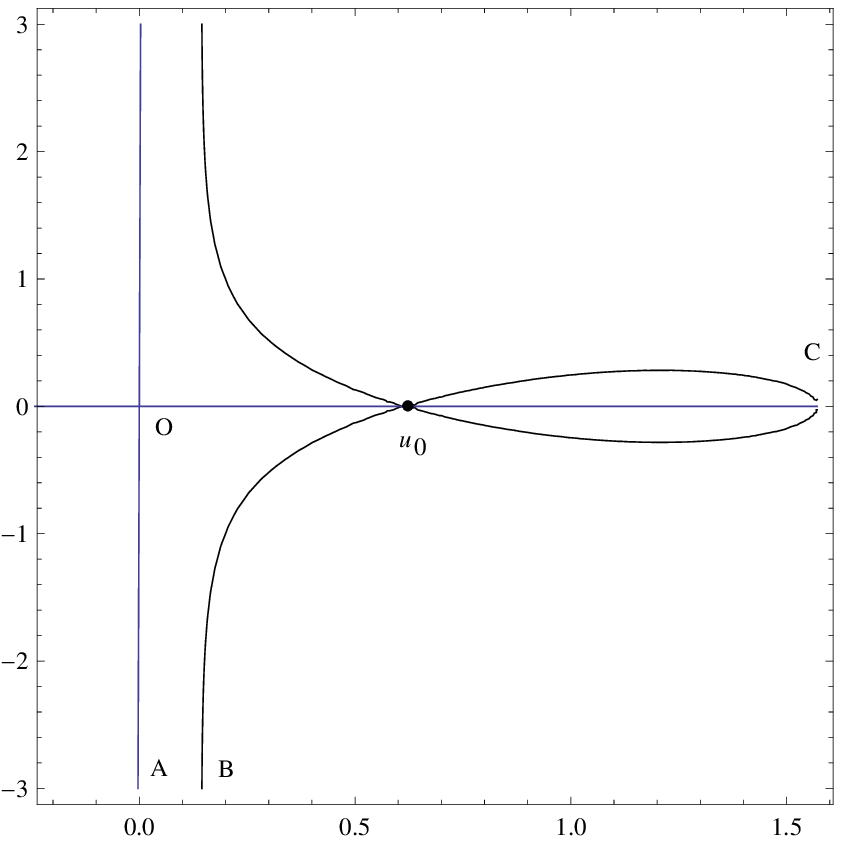}
	\qquad
	{\tiny($b$)}\ \includegraphics[width=0.35\textwidth]{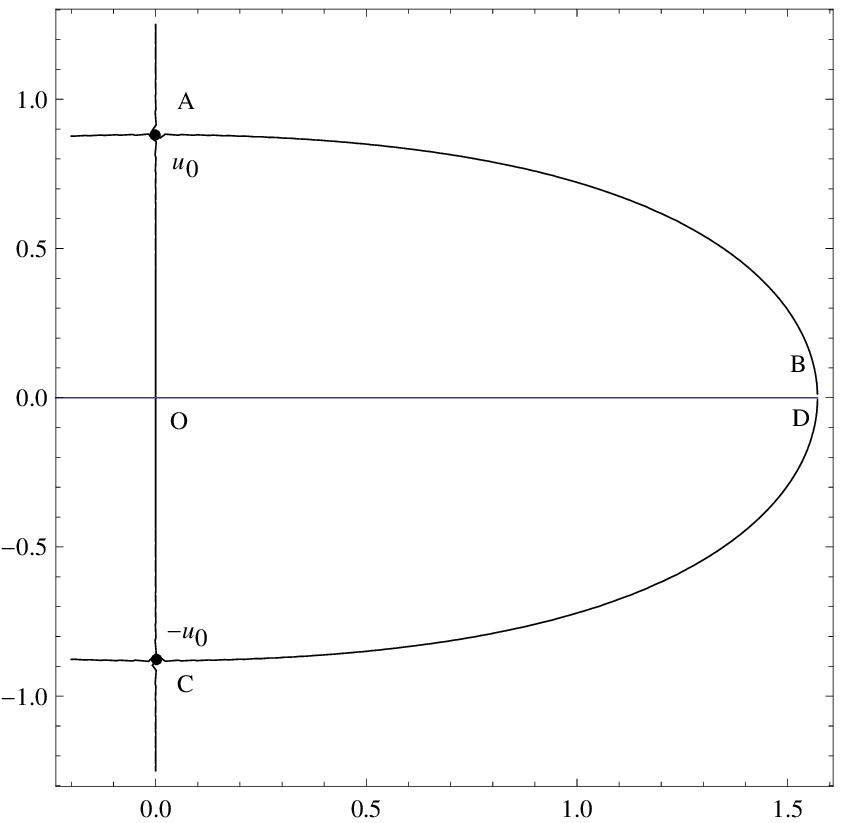} 
\caption{\small{Integration paths through the saddle $u_0$: (a) $a>1$ and (b) $a<1$.}}
\end{center}
\end{figure}

\vspace{0.3cm}

\noindent{\bf 3.1}\ \ {\bf The case $a>1$}
\vspace{0.3cm}

\noindent The steepest descent path for the exponential factor $\exp\{x\psi(u)\}$ is shown in Fig.~2(a). This consists of the path $OA$ along the negative imaginary axis followed by the path $BC$ through the saddle $u_0=\arctan \sqrt{a-1}$ (this last path is the same as that in Fig.~1(a)). The integral with the exponential factor $\exp\{-x\psi(u)\}$ is taken along the mirror image of these paths and yields a conjugate expansion.
Then
\[\frac{1}{\pi i} \int_0^{\pi/2} e^{x\psi(u)}du=I_{OA}+I_{BC},\] 
where
\[I_{OA}=\frac{1}{\pi i}\int_0^{-\infty i}e^{x\psi(u)}du,\qquad
I_{BC}=\frac{1}{\pi i}\int_{BC}e^{x\psi(u)}du\]
and the contribution from the horizontal connecting path $AB$ vanishes as $\Im (u)\to-\infty$.

To deal with the integral $I_{OA}$, we put $u=-iy$ to obtain
\[I_{OA}=-\frac{1}{\pi}\int_0^\infty e^{-x(ay-\tanh y)}dy=-\frac{1}{\pi}\int_0^\infty e^{-x(a-1) w} \frac{dy}{dw}\,dw,\]
where we have introduced the new variable $w$ by
\[(a-1)w=ay-\tanh y=(a-1)y+\frac{1}{3}y^3-\frac{2}{15}y^5+\frac{17}{315}y^7+\cdots\ .\]
Inversion of this series followed by differentiation yields (no odd powers of $w$ are present)
\bee\label{e31a}
\frac{dy}{dw}=\sum_{k\geq0} C_k(a) w^{2k},
\ee
where the first few coefficients $C_k(a)$ are:
\[C_0(a)=1,\quad C_1(a)=\frac{1}{1-a},\quad C_2(a)=\frac{3+2a}{3(1-a)^2},\quad
 C_3(a)=\frac{(45 + 78 a + 17 a^2)}{45(1-a)^3},\]
 \[\quad C_4(a)= \frac{(315 + 972 a + 576 a^2 + 62 a^3)}{315(1-a)^4},\quad
C_5(a)=\frac{(14175 + 66060 a + 71982 a^2 + 21576 a^3 + 1382 a^4)}{14175(1-a)^5},\]
\bee\label{e31b}
C_6(a)=\frac{(467775 + 3001590 a + 5063616 a^2 + 2842542 a^3 + 514533 a^4 + 
  21844 a^5)}{467775(1-a)^6}.
\ee
Then
\[
I_{OA}\sim -\frac{1}{\pi}\sum_{k=0}^\infty C_k(a)\int_0^\infty e^{-\la w}w^{2k}dw
= -\frac{1}{\pi\la}\sum_{k=0}^\infty \frac{C_k(a)(2k)!}{\la^{2k}},\]
where $\la:=x(a-1)=\nu-x$.

The contribution $I_{BC}$ can be obtained from (\ref{e22}) in the form
\[I_{BC}\sim \frac{e^{i(x\phi(u_0)+\pi/4)}}{i\sqrt{\pi ax}}\,(a-1)^{-1/4}\sum_{k=0}^\infty \frac{i^k \g(k+\fs)}{\g(\fs)}\,\frac{A_{2k}(a)}{x^k},\]
where the coefficients $A_{2k}(a)$ are listed in (\ref{e23b}). Taking into account the conjugate expansion resulting from the exponential factor $\exp\{-x\psi(u)\}$ taken along the mirror image of the paths in Fig.~2(a), we finally end up with the following expansion:
\begin{theorem}$\!\!\!.$\ \ The asymptotic expansion of $h_\nu(x)$ for $x\to+\infty$ with $\nu=ax$, $a>1$ is
\[h_\nu(x)\sim \frac{2(a-1)^{-1/4}}{\sqrt{\pi ax}}\bl\{\sin \Phi(x)\sum_{k=0}^\infty \frac{(-)^kA_{4k}(a) (\fs)_{2k}}{x^{2k}} +\frac{\cos \Phi(x)}{2x}
\sum_{k=0}^\infty \frac{(-)^k A_{4k+2}(a) (\f{3}{2})_{2k}}{x^{2k}}\br\}\]
\bee\label{e31}
-\frac{2}{\pi x(a-1)} \sum_{k=0}^\infty \frac{(2k)!}{x^{2k}}\,\frac{C_k(a)}{(a-1)^{2k}}
\ee
where $\Phi(x)$ is defined in (\ref{e23a}). The coefficients $A_{2k}(a)$ and $C_k(a)$ are defined in (\ref{e23b}) and (\ref{e31b}).
\end{theorem}

\vspace{0.3cm}

\noindent{\bf 3.2}\ \ {\bf The case $a<1$}
\vspace{0.3cm}

\noindent The steepest descent path for the exponential factor $\exp \{x\psi(u)\}$ in the case $a<1$ is shown in Fig.~2(b). This consists of the positive imaginary axis $OA$ between $[0,|u_0|]$, where $|u_0|=\mbox{arctanh} \sqrt{1-a}$, 
together with the steepest path $AB$ emanating from the saddle and terminating at the endpoint $u=\fs\pi$. The integral with the exponential factor $\exp\{-x\psi(u)\}$ is taken over the image of this path in the lower half plane to produce the conjugate result.

Then 
\[\frac{1}{\pi i}\int_0^{\pi/2}e^{x\psi(u)}du =\frac{1}{\pi i}\bl\{\int_{OA}+\int_{AB}\br\}e^{x\psi(u)}\,du,\]
where
\[I_{OA}=\frac{1}{\pi i}\int_0^{i|u_0|} e^{x\psi(u)} du=\frac{1}{\pi}\int_0^{|u_0|}e^{-x(\tanh y-ay)}dy.\]
We now introduce the new variable $w$ in a similar manner to that in \S 3.1 to find
%\[(1-a)w=\tanh y-ay=(1-a)y-\frac{1}{3}y^3+\frac{2}{15}y^5-\frac{17}{315}y^7+\cdots\]
%to find upon inversion and differentiation the expansion (no odd powers of $w$)
%\[\frac{dy}{dw}=\sum_{k\geq 0} {\hat e}_k w^{2k},\]
%where
%\[{\hat e}_0=1,\quad {\hat e}_1=\frac{1}{1-a},\quad {\hat e}_2=\frac{3+2a}{3(1-a)^2},\quad
% {\hat e}_3=\frac{(45 + 78 a + 17 a^2)}{45(1-a)^3},\]
% \[\quad {\hat e}_4= \frac{(315 + 972 a + 576 a^2 + 62 a^3)}{315(1-a)^4},\quad
%{\hat e}_5=\frac{(14175 + 66060 a + 71982 a^2 + 21576 a^3 + 1382 a^4)}{14175(1-a)^5},\]
%\[{\hat e}_6=\frac{(467775 + 3001590 a + 5063616 a^2 + 2842542 a^3 + 514533 a^4 + 
%  21844 a^5)}{467775(1-a)^6}.\]
\[I_{OA}=\frac{1}{\pi}\int_0^{w_0}e^{-\la w} \frac{dy}{dw}\,dw=\frac{1}{\pi}\sum_{k=0}^\infty C_k(a) \int_0^{w_0} e^{-\la w} w^{2k}dw,\]  
where now $\la=x(1-a)$ and $w_0:=(\tanh |u_0|-a|u_0|)/(1-a)$. Evaluation of the integrals in terms of the (lower) incomplete gamma function $\gamma(\alpha,z)=\int_0^z e^{-t}t^{\alpha-1}dt$ ($\Re (\alpha)>0$) then yields
\bee\label{e32}
I_{OA}=\frac{1}{\pi\la}\sum_{k=0}^\infty \frac{C_k(a)}{\la^{2k}}\,\gamma(2k+1,x\Psi(a)),
\ee
with $\la w_0=x\Psi(a)$, where $\Psi(a)$ is defined in (\ref{e23a}). A similar argument for the integral involving factor $\exp \{-x\psi(u)\}$ taken along the path $[0,-i|u_0|]$ shows that
\bee\label{e32a}
\frac{1}{\pi i}\int_0^{-i|u_0|}e^{-x\psi(u)}du=-I_{OA}.
\ee

The determination of the expansion of $I_{AB}$ is similar to that described in \S 2.2, with one important difference.
Because we have the difference between the contributions from the steepest descent paths $AB$ (with $\exp\{x\psi(u)\}$) and $CD$ (with $\exp\{-x\psi(u)\}$), we find
\[I_{AB}-I_{CD}=\frac{1}{\pi i}\int_{AB}e^{x\psi(u)}du-\frac{1}{\pi i}\int_{CD}e^{-x\psi(u)}du\]
\[=\frac{e^{-x\Psi(a)}}{\pi i} \int_0^\infty e^{-xw^2} \bl(\frac{du}{dw}-{\overline\frac{du}{dw}}\br) dw,\]
where the expansion of $du/dw$ given in (\ref{e23c}) now has to include the odd-order powers of $w$ since integration is along only half of the steepest descent path.
Thus we have
\[\frac{du}{dw}=\frac{(1-a)^{-1/4}}{\sqrt{a}}\bl\{\sum_{k\geq0} {\hat A}_{2k}(a)w^{2k}+i\sum_{k\geq0} (-)^k {\hat A}_{2k+1}(a)w^{2k+1}\br\},\]
where the first few odd-order coefficients are:
\[{\hat A}_1(a)=\frac{3a-2}{3a(1-a)},\quad {\hat A}_3(a)=\frac{4(9a-4)}{135a^2(1-a)^{5/2}},\quad {\hat A}_5(a)=
\frac{32-120a+144a^2+189a^3}{2835a^3(1-a)^4},\]
\[{\hat A}_7(a)=\frac{8(16-84a+180a^2-207a^3+270a^4)}{25515a^4(1-a)^{11/2}},\]
\bee\label{e33b}
{\hat A}_9(a)=\frac{(-35968 + 242784 a - 692064 a^2 + 1077948 a^3 - 1020600 a^4 + 
 1403325 a^5 + 400950 a^6)}{37889775a^5(1-a)^7}~.
\ee
 
Then
\[I_{AB}-I_{CD}\sim \frac{2(1-a)^{-1/4}}{\pi\sqrt{a}}\,e^{-x\Psi(a)} \sum_{k=0}^\infty (-)^k {\hat A}_{2k+1}(a) \int_0^\infty e^{-xw^2} w^{2k+1} dw\]
\bee\label{e33}
=\frac{(1-a)^{-1/4}}{\pi\sqrt{ax}}e^{-x\Psi(a)} \sum_{k=0}^\infty \frac{(-)^k k! {\hat A}_{2k+1}(a)}{x^k}.
\ee
Collecting together the results in (\ref{e32}), (\ref{e32a}) and (\ref{e33}), we finally obtain:
\begin{theorem}$\!\!\!.$\ \ The asymptotic expansion of $h_\nu(x)$ for $x\to+\infty$ with $\nu=ax$, $a<1$ is
\[h_\nu(x)\sim \frac{(1-a)^{-1/4}}{\pi\sqrt{ax}}\,e^{-x\Psi(a)} \sum_{k=0}^\infty \frac{(-)^k k! {\hat A}_{2k+1}(a)}{x^k}\hspace{4cm}\]
\bee\label{e34}
\hspace{4cm}+\frac{2}{\pi x(1-a)} \sum_{k=0}^\infty \frac{C_k(a)}{x^{2k}}\,\frac{\gamma(2k+1,x\Psi(a))}{(1-a)^{2k}},
\ee
where $\Psi(a)$ is defined in (\ref{e24a}). The coefficients ${\hat A}_k(a)$ and $C_k(a)$ are defined in (\ref{e23b}), (\ref{e33b}) and (\ref{e31b}).
\end{theorem}
\vspace{0.3cm}

\noindent{\bf 3.3}\ \ {\bf The case $a=1$}
\vspace{0.3cm}

\noindent The steepest descent paths when $a=1$ are shown in Fig.~1(c). From (\ref{e26}), we obtain:
\begin{theorem}$\!\!\!.$\ \ The asymptotic expansion of $h_\nu(\nu)$ for $\nu\to+\infty$ is
\bee\label{e35}
h_\nu(\nu)=\frac{1}{i}(I_{OA}-{\overline I}_{OA})\sim\frac{2}{3\pi \nu^{1/3}} \sum_{k=0}^\infty \frac{\Im (B_k)}{\nu^{2k/3}}\,\g(\f{2}{3}k+\f{1}{3}).
\ee
Substitution of the values of the coefficients $B_k$ from (\ref{e2dk}) followed by some routine algebra then shows that
\[h_\nu(\nu)\sim \frac{3^{-2/3}}{\pi x^{1/3}}\bl\{\g(\f{1}{3})+\frac{9\cdot 3^{1/3}}{35 \nu^{4/3}} \g(\f{5}{3})+\frac{2}{35\nu^2} \g(\f{7}{3})-\frac{3258\cdot 3^{1/3}}{284375 \nu^{10/3}} \g(\f{11}{3})-\frac{49711}{11790625\nu^4} \g(\f{13}{3})+\cdots\br\}\]
\bee\label{e36}
-\frac{4}{5\pi \nu}
\bl\{1-\frac{216\cdot 2!}{13475 \nu^2}+\frac{565056\cdot 4!}{372246875\nu^4}-\frac{102278311479936\cdot 6!}{712401183634765625 \nu^6}+\cdots\br\}
\ee
valid as $\nu\to+\infty$.
\end{theorem}
\vspace{0.6cm}

\begin{center}
{\bf 4.\ The asymptotic expansions of $k_\nu(x)$ and $h_\nu(x)$ when $x\to-\infty$, $\nu=a|x|$}
\end{center}
\setcounter{section}{4}
\setcounter{equation}{0}
\renewcommand{\theequation}{\arabic{section}.\arabic{equation}} 
To deal with the functions with negative argument, we replace $x$ by $-x$, where $x>0$, and again set $\nu=ax$, $a>0$ to find from (\ref{e11}) and (\ref{e12})
\[k_\nu(-x)=\frac{2}{\pi}\int_0^{\pi/2} \cos (x(\tan u+au))\,du=\frac{1}{\pi}\int_0^{\pi/2}\{e^{x\psi(u)}+e^{-x\psi(u)}\}du,\]
\[h_\nu(-x)=-\frac{2}{\pi }\int_0^{\pi/2} \sin (x(\tan u+au))\,du=-\frac{1}{\pi i}\int_0^{\pi/2}\{e^{x\psi(u)}-e^{-x\psi(u)}\}du,\]
where now 
\[\psi(u)=i(\tan u+au).\]
Saddle points occur at $\psi'(u)=i(\mbox{sec}^2u+a)=0$; that is, in the domain of interest at the points
\[u_0^\pm=\fs\pi\pm i\beta,\quad \beta=\mbox{arctanh}\, 1/\sqrt{1+a}.\]

Since $\psi''(u_0^\pm)=2a\sqrt{1+a}>0$, the steepest descent path through the saddles is locally parallel to the imaginary axis. In fact, it is easy to see that, for $y>0$, $\Im \psi(\fs\pi+iy)=\Im \psi(\fs\pi+i\beta)=\fs\pi a$, so that the steepest descent path through the saddle in the upper half-plane is the vertical line $\Re (u)=\fs\pi$, $\Im (u)\geq 0$ for the integral with $\exp \{x\psi(u)\}$; similarly, for the integral with $\exp\{-x\psi(u)\}$ the steepest path is $\Re (u)=\fs\pi$, $\Im (u)\leq 0$. The integration path $[0, \fs\pi]$ is then deformed into the path $OABC$ for the integral involving $\exp\{x\psi(u)\}$ and $OA'B'C$ for the integral involving $\exp\{-x\psi(u)\}$; see Fig.~3. The contributions from the ends $AB$ and $A'B'$ vanish as $\Im (u)\to\pm\infty$.

\begin{figure}[th]
	\begin{center}\ \includegraphics[width=0.25\textwidth]{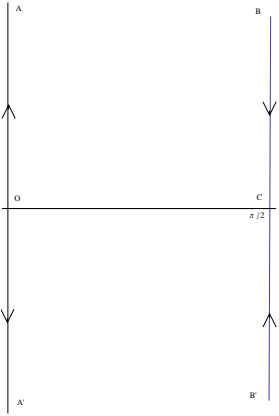}	
\caption{\small{Integration paths for the case $x\to -\infty$.}}
\end{center}
\end{figure}

We deal with the integral containing the factor $\exp\{x\psi(u)\}$. Introducing the new variable
\[-w^2=\psi(u)-\psi(u_0^+)=a\sqrt{1+a} (u-u_0^+)^2+ \frac{i}{3}a(2+3a) (u-u_0^+)^3+\cdots,\]%-\frac{1}{3}a\sqrt{1+a}(1+3a)(y-\beta)^4+\cdots\, ,\]
we find upon substitution of $u=\fs\pi+iy$, $y\geq0$ (so that $u-u_0^+=i(y-\beta$)) and inversion that (compare (\ref{e23c}))
\[\frac{dy}{dw}\stackrel{e}{=}\frac{(1+a)^{-1/4}}{\sqrt{a}} \sum_{k\geq0} {\hat A}_{2k}(-a) w^{2k},\]
where the coefficients ${\hat A}_{2k}(-a)$ can be obtained from (\ref{e23d}). Then the contribution from the path $BC$ is
\[I_{BC}=\frac{1}{\pi}\int_{\pi/2+\infty i}^{\pi/2}e^{x\psi(u)}du=\frac{1}{\pi i} \int_0^\infty e^{x\psi(\frac{1}{2}\pi+iy)}dy
=\frac{e^{x\psi(u_o^+)}}{\pi i}\int_{-\infty}^\infty e^{-xw^2} \frac{dy}{dw}\,dw\]
\bee\label{e41}
\sim\frac{(1+a)^{-1/4}}{i\sqrt{\pi ax}}\,e^{x\psi(u_0^+)} \sum_{k=0}^\infty \frac{\g(k+\fs)}{\g(\fs)}\,\frac{{\hat A}_{2k}(-a)}{x^k},
\ee
where
\bee\label{e41a}
x\psi(u_0^+)=\fs\pi i\nu-x\Omega(a),\qquad \Omega(a):=\sqrt{1+a}+a \,\mbox{arctanh}\,1/\sqrt{1+a}.
\ee
The contribution from the path $B'C$ is given by the conjugate expansion since
\bee\label{e42}
I_{B'C}=\frac{1}{\pi}\int_{\pi/2-\infty i}^{\pi/2} e^{-x\psi(u)}du=-\frac{1}{\pi i}\int_0^\infty e^{-x\psi(\frac{1}{2}\pi-iy)}dy={\overline I}_{BC}.
\ee

We now consider the contributions from the paths $OA$ and $OA'$ along the imaginary axis. We have
\[I_{OA}=\frac{1}{\pi}\int_0^{\infty i} e^{x\psi(u)}du=\frac{i}{\pi}\int_0^\infty e^{-x(ay+\tanh y)}dy
%=\frac{1}{\pi}\int_0^{-\infty i} e^{-x\psi(u)}du=-I_{OA'}\]
=\frac{i}{\pi}\int_0^\infty e^{-\la w}\,\frac{dy}{dw}\,dw,\qquad \la:=x(1+a),\]
where, from (\ref{e31a}) and (\ref{e31b}),
\[\frac{dy}{dw}=\sum_{k\geq 0}C_k(-a) w^{2k}.\]
Hence we obtain the expansion
\bee\label{e43}
I_{OA}=-I_{OA'}\sim\frac{i}{\pi\la}\sum_{k=0}^\infty \frac{(2k)! C_k(-a)}{\la^{2k}},
\ee
since it is easily seen that $I_{OA'}=-I_{OA}$.

Collecting together the results in (\ref{e41}) -- (\ref{e43}), we finally obtain the expansions:
\begin{theorem}$\!\!\!.$\ \ The asymptotic expansions of $k_\nu(-x)$ and $h_\nu(-x)$ for $x\to+\infty$ with $\nu=ax$, $a>0$ are
\begin{eqnarray}k_\nu(-x)&=&I_{OA}+I_{OA'}+I_{BC}+I_{B'C}\nonumber\\
&\sim& \frac{2(1+a)^{-1/4}}{\sqrt{\pi ax}} \, e^{-x\Omega(a)}\sin \fs\pi\nu \sum_{k=0}^\infty \frac{\g(k+\fs)}{\g(\fs)}\,\frac{{\hat A}_{2k}(a)}{x^k}\label{e44}
\end{eqnarray}
and
\begin{eqnarray}
h_\nu(-x)&=&i\{I_{OA}-I_{OA'}+I_{BC}-I_{B'C}\}\nonumber\\
&\sim& \frac{2(1+a)^{-1/4}}{\sqrt{\pi ax}} \, e^{-x\Omega(a)}\cos \fs\pi\nu \sum_{k=0}^\infty \frac{\g(k+\fs)}{\g(\fs)}\,\frac{{\hat A}_{2k}(a)}{x^k}\nonumber\\
&&\hspace{2cm} -\frac{2}{\pi x(1+a)}\sum_{k=0}^\infty \frac{(2k)!\, C_k(-a)}{x^k(1+a)^k}.\label{e45}
\end{eqnarray}
\end{theorem}
It may be observed that the first asymptotic series in (\ref{e45}) is exponentially small (and vanishes when $\nu=1, 3, 5, \ldots$) and that second series is the dominant
expansion of O($x^{-1})$. 
\vspace{0.6cm}

\begin{center}
{\bf 5.\ Numerical results}
\end{center}
\setcounter{section}{5}
\setcounter{equation}{0}
\renewcommand{\theequation}{\arabic{section}.\arabic{equation}} 
We present numerical results to confirm the accuracy of the various expansions obtained. In Table 1 values of the Bateman function $k_\nu(x)$ are shown for two values of $x$ and values of the parameter $a$ corresponding to the three cases $a>1$, $a<1$ and $a=1$ considered in Section 2. The expansions employed are those given in (\ref{e23}), (\ref{e24}) and (\ref{e25}) using the number of coefficients displayed. The exact value of $k_\nu(x)$ was obtained from the hypergeometric representation in (\ref{e13}) and the absolute relative error in the different expansions is also indicated.  
\begin{table}[h]
\caption{\footnotesize{Values of the Bateman function $k_\nu(x)$ with $\nu=ax$.}}
\begin{center}
\begin{tabular}{|c|c|c|}
\hline
&&\\[-0.3cm]
%\mcol{1}{|c|}{} & \mcol{2}{c|}{$n=20$} & \mcol{2}{c|}{$n=40$}\\
\mcol{1}{|c|}{} & \mcol{1}{c|}{$a=2,\ x=20$} & \mcol{1}{c|}{$a=2,\ x=60$}\\
\hline
&&\\[-0.3cm]
$k_\nu(x)$ & $-6.5410626744\times 10^{-2}$ & $-4.6850153400\times 10^{-2}$ \\
Asymptotic & $-6.5410542459\times 10^{-2}$ & $-4.6850153626\times 10^{-2}$ \\
Error      & $1.289\times 10^{-6}$ & $4.830\times 10^{-9}$\\
\hline
\mcol{1}{|c|}{} & \mcol{1}{c|}{$a=0.50,\ x=20$} & \mcol{1}{c|}{$a=0.50,\ x=60$}\\
\hline
&&\\[-0.3cm]
$k_\nu(x)$ & $+1.0048261319\times 10^{-3}$ & $+1.3875830492\times 10^{-8}$ \\
Asymptotic & $+1.0048449917\times 10^{-3}$ & $+1.3875831812\times 10^{-8}$ \\
Error      & $1.877\times 10^{-5}$ & $8.510\times 10^{-8}$\\
\hline
\mcol{1}{|c|}{} & \mcol{1}{c|}{$a=1,\ x=20$} & \mcol{1}{c|}{$a=1,\ x=60$}\\
\hline
&&\\[-0.3cm]
$k_\nu(x)$ & $+2.6100825169\times 10^{-1}$ & $+1.8127952802\times 10^{-1}$ \\
Asymptotic & $+2.6100825257\times 10^{-1}$ & $+1.8127952803\times 10^{-1}$ \\
Error      & $3.383\times 10^{-9}$ & $9.272\times 10^{-12}$\\
\hline
\end{tabular}
\end{center}
\end{table}
\begin{table}[h]
\caption{\footnotesize{Values of the Havelock function $h_\nu(x)$ with $\nu=ax$.}}
\begin{center}
\begin{tabular}{|c|c|c|}
\hline
&&\\[-0.3cm]
%\mcol{1}{|c|}{} & \mcol{2}{c|}{$n=20$} & \mcol{2}{c|}{$n=40$}\\
\mcol{1}{|c|}{} & \mcol{1}{c|}{$a=2,\ x=20$} & \mcol{1}{c|}{$a=2,\ x=60$}\\
\hline
&&\\[-0.3cm]
$h_\nu(x)$ & $+1.3427850086\times 10^{-1}$ & $-1.0233792538\times 10^{-1}$ \\
Asymptotic & $+1.3427777585\times 10^{-1}$ & $-1.0233792529\times 10^{-1}$ \\
Error      & $5.399\times 10^{-6}$ & $8.881\times 10^{-10}$\\
\hline
\mcol{1}{|c|}{} & \mcol{1}{c|}{$a=0.50,\ x=20$} & \mcol{1}{c|}{$a=0.50,\ x=60$}\\
\hline
&&\\[-0.3cm]
$h_\nu(x)$ & $+6.7009789898\times 10^{-2}$ & $+2.1318746869\times 10^{-2}$ \\
Asymptotic & $+6.6562158637\times 10^{-2}$ & $+2.131875597\times 10^{-2}$ \\
Error      & $1.877\times 10^{-5}$ & $4.271\times 10^{-7}$\\
\hline
\mcol{1}{|c|}{} & \mcol{1}{c|}{$a=1,\ x=20$} & \mcol{1}{c|}{$a=1,\ x=60$}\\
\hline
&&\\[-0.3cm]
$h_\nu(x)$ & $+1.3865680166\times 10^{-1}$ & $+1.0052911205\times 10^{-1}$ \\
Asymptotic & $+1.3865681475\times 10^{-1}$ & $+1.0052911217\times 10^{-1}$ \\
Error      & $9.439\times 10^{-8}$ & $1.159\times 10^{-9}$\\
\hline
\end{tabular}
\end{center}
\end{table}

Table 2 shows the same values of $a$ and $x$ for the Havelock function $h_\nu(x)$. The asymptotic values were obtained from (\ref{e31}), (\ref{e34}) and (\ref{e36}). The numerical evaluation of $h_\nu(x)$ for large $x$ cannot be evaluated directly from the integral (\ref{e12}) on account of the fact that the integrand is highly oscillatory as $u\to\fs\pi$. Numerical values of $h_\nu(x)$ were obtained by integrating along the steepest descent paths shown in Figs. 1 and 2, thereby eliminating the rapid oscillations (we omit details of these calculations).

\begin{table}[h]
\caption{\footnotesize{Values of the functions $k_\nu(-x)$ and $h_\nu(-x)$ with $\nu=ax$.}}
\begin{center}
\begin{tabular}{|c|c|c|}
\hline
&&\\[-0.3cm]
%\mcol{1}{|c|}{} & \mcol{2}{c|}{$n=20$} & \mcol{2}{c|}{$n=40$}\\
\mcol{1}{|c|}{} & \mcol{1}{c|}{$a=0.25,\ x=10$} & \mcol{1}{c|}{$a=0.75,\ x=20$}\\
\hline
&&\\[-0.3cm]
$k_\nu(-x)$& $-1.9280268893\times 10^{-7}$ & $-2.1016307406\times 10^{-19}$ \\
Asymptotic & $-1.9278470013\times 10^{-7}$ & $-2.1016310895\times 10^{-19}$ \\
Error      & $9.330\times 10^{-5}$ & $1.125\times 10^{-7}$\\
\hline
\mcol{1}{|c|}{} & \mcol{1}{c|}{$a=1.50,\ x=10$} & \mcol{1}{c|}{$a=1.75,\ x=30$}\\
\hline
&&\\[-0.3cm]
$k_\nu(-x)$& $-4.4526301254\times 10^{-13}$ & $-1.4514839956\times 10^{-26}$ \\
Asymptotic & $-4.4526319222\times 10^{-13}$ & $-1.4514840241\times 10^{-26}$ \\
Error      & $6.680\times 10^{-3}$ & $1.958\times 10^{-8}$\\
\hline
\mcol{1}{|c|}{} & \mcol{1}{c|}{$a=0.25,\ x=10$} & \mcol{1}{c|}{$a=0.75,\ x=15$}\\
\hline
&&\\[-0.3cm]
$h_\nu(-x)$& $-5.1481358274\times 10^{-2}$ & $-2.4292604794\times 10^{-2}$ \\
Asymptotic & $-5.1481463710\times 10^{-2}$ & $-2.4292604793\times 10^{-2}$ \\
Error      & $2.048\times 10^{-6}$ & $1.445\times 10^{-11}$\\
\hline
\mcol{1}{|c|}{} & \mcol{1}{c|}{$a=1.50,\ x=10$} & \mcol{1}{c|}{$a=1.75,\ x=25$}\\
\hline
&&\\[-0.3cm]
$h_\nu(-x)$& $-2.5497382200\times 10^{-2}$ & $-1.5439800067\times 10^{-2}$ \\
Asymptotic & $-2.5497382199\times 10^{-2}$ & $-1.5439800067\times 10^{-2}$ \\
Error      & $2.817\times 10^{-11}$ & $7.313\times 10^{-15}$\\
\hline
\end{tabular}
\end{center}
\end{table}

Finally, in Table 3 we show the values of $k_\nu(-x)$ and $h_\nu(-x)$ and their asymptotic representations in (\ref{e44}) and (\ref{e45}) when $\nu=ax$. The exact evaluation of $k_\nu(-x)$ was obtained from (\ref{e14}) and that of $h_\nu(-x)$ was determined by integration along the paths shown in Fig.~3.

\vspace{0.6cm}

\begin{center}
{\bf 6.\ Concluding remarks}
\end{center}
\setcounter{section}{6}
\setcounter{equation}{0}
\renewcommand{\theequation}{\arabic{section}.\arabic{equation}} 
Asymptotic expansions have been derived for the Bateman and Havelock functions $k_\nu(x)$ and $h_\nu(x)$ for $x\to\pm\infty$ when the order $\nu=a|x|$, with the parameter $a>0$, by application of the method of steepest descents.
In the case $x\to+\infty$, three different cases arise according as $a<1$, $a=1$ and $a>1$. In the case $x\to-\infty$. a single expansion holds for $a>0$. The coefficients in these expansions are obtained using an inversion procedure (essentially Lagrange inversion). It is worth mentioning that, although the complexity of the coefficients rapidly increases with increasing order $k$, in the case of specific numerical values for $a=\nu/x$ 
it is possible to generate many more coefficients than those displayed here.
Numerical results employing the coefficients listed in the paper show excellent agreement with the numerically computed values of $k_\nu(\pm x)$ and $h_\nu(\pm x)$.

In the case of the Bateman function, a uniform approximation valid when $a\simeq 1$ and $x>0$ can be obtained in terms of the Airy function Ai from (\ref{e13}) and \cite[p.~412, Ex.~7.3]{Olv}. If we define the quantity $\zeta$ by
\begin{eqnarray*}
\f{2}{3}\zeta^{3/2}\!\!&=&\!\!\sqrt{1-a}-a\,\mbox{arctanh} \sqrt{1-a}\qquad ( 0<a\leq 1,\ \zeta\geq0)\\
\f{2}{3}(-\zeta)^{3/2}\!\!&=&\!\!a\,\arctan \sqrt{a-1}-\sqrt{a-1}\qquad (a\geq 1,\ \zeta\leq0),
\end{eqnarray*}
we obtain the approximation
\bee\label{e61}
k_\nu(x)\sim \frac{2}{a^{1/2} x^{1/3}} \bl(\frac{\zeta}{1-a}\br)^{1/4} \mbox{Ai} (x^{2/3}\zeta)\qquad (x\to+\infty),
\ee
when $\nu=ax$ with $a\simeq 1$. Use of the asymptotic behaviour of Ai($z)$ given by
\[\mbox{Ai} (z)\sim \frac{e^{-z}}{2\sqrt{\pi}z^{1/4}},\qquad \mbox{Ai} (-z)\sim \frac{1}{\sqrt{\pi} z^{1/4}}\,\cos \bl(\frac{2}{3}z^{3/2}-\frac{\pi}{4}\br)\qquad (z\to+\infty),\]
combined with Ai$(0)=\g(\f{1}{3})/(2\cdot 3^{1/6}\pi)$ and the fact that $\zeta/(1-a)\to1$ as $a\to 1^-$, shows that (\ref{e61}) reduces to the leading terms of the expansions given Theorems 1, 2 and 3 when $a>1$, $a<1$ and $a=1$.

\vspace{0.6cm}

\begin{center}
{\bf Appendix A:  Asymptotic expansion of $h_\nu(\pm x)$ for $x\to+\infty$ and finite $\nu>0$}
\end{center}
\setcounter{section}{1}
\setcounter{equation}{0}
\renewcommand{\theequation}{\Alph{section}.\arabic{equation}}
We deform the integration path $[0,\fs\pi]$ in (\ref{e12}) into the path depicted in Fig.~2(b), where with $x\to+\infty$, $\nu$ fixed the parameter $a=\nu/x\to 0$. The saddle point $u_0=i\,\mbox{arctanh}\,\sqrt{1-a}\sim i\log\,(2/\sqrt{a})$ as $a\to 0$. Then the Havelock function is
\bee\label{a1}
h_\nu(x)=\frac{1}{\pi i} \int_0^{\pi/2}\{e^{i(x\tan u-\nu u)}+e^{-i(x\tan u-\nu u)}\}\,du.
\ee

Consider the integral involving the first exponential in (\ref{a1}) taken round the path $OAB$ in Fig.~2(b):
\[\frac{1}{\pi i}\int_0^{\pi/2} e^{i(x\tan u-\nu u)}du=I_{OA}+I_{AB},\]
where, from (\ref{e33}), $I_{AB}=O(e^{-x\Psi(a)})=O(e^{-x})$ as $a\to 0$. The integral along the positive imaginary axis becomes
\[I_{OA}=\frac{1}{\pi}\int_0^{|u_0|} e^{-x\tanh s+\nu s} ds=\frac{1}{\pi}\int_0^{w_0} e^{-xw+\nu {\footnotesize\mbox{arctanh}}\,w}\,\frac{dw}{1-w^2},\]
where we have made the substitution $s=\tanh w$ and $w_0=\sqrt{1-a}<1$. Using the expansion
\[\frac{e^{\pm\nu {\footnotesize\mbox{arctanh}}\,w}}{1-w^2}=\sum_{k=0}^\infty \frac{(\pm)^kc_k(\nu)}{k!}\,w^k\qquad (|w|<1),\]
where the first few coefficients $c_k(\nu)$ are
\begin{eqnarray*}
c_0(\nu)&=&1,\quad c_1(\nu)=\nu,\quad c_2(\nu)=2+\nu^2,\quad c_3(\nu)=\nu(8+\nu^2),\\
c_4(\nu)&=&24+20\nu^2+\nu^4,\quad
c_5(\nu)=\nu(184+40\nu^2+\nu^4),\\
c_6(\nu)&=&720+784\nu^2+70\nu^4+\nu^6, \quad
c_7(\nu)=\nu(8448+2464\nu+112\nu^4+\nu^6),\\
c_8(\nu)&=&40320+52352\nu^2+6384\nu^4+168\nu^6+\nu^8,
\end{eqnarray*}
we find
\[I_{OA}=\frac{1}{\pi}\sum_{k=0}^\infty \frac{c_k(\nu)}{k!} \int_0^{w_0} e^{-xw} w^kdw\sim\frac{1}{\pi x} \sum_{k=0}^\infty \frac{c_k(\nu)}{x^k}\qquad (x\to+\infty).\]
In the final step the upper limit $w_0$ has been replaced by $\infty$.

The integral involving the second exponential in (\ref{a1}) taken along the mirror image in the lower-half plane of the path $OAB$ in Fig.~2(b)  yields the same result. Hence we obtain the expansion
\bee\label{a2}
h_\nu(x)\sim \frac{2}{\pi x} \sum_{k=0}^\infty \frac{c_k(\nu)}{x^k} \qquad (x\to+\infty,\ \nu\ \mbox{fixed}).
\ee
The asymptotic expansion of $h_\nu(-x)$ for $x\to+\infty$ with $\nu$ fixed follows the same procedure to yield
\bee\label{a3}
h_\nu(-x)\sim \frac{2}{\pi x} \sum_{k=0}^\infty \frac{(-)^{k-1}c_k(\nu)}{x^k} \qquad (x\to+\infty,\ \nu\ \mbox{fixed}),
\ee
where in both (\ref{a2}) and (\ref{a3}) we have neglected an exponentially small contribution of $O(e^{-x})$. 
\vspace{0.6cm}

\begin{center}
{\bf Appendix B:  The steepest descent paths in Fig.~1(a) when $a>1$}
\end{center}
\setcounter{section}{2}
\setcounter{equation}{0}
\renewcommand{\theequation}{\Alph{section}.\arabic{equation}}
The steepest descent paths through the saddle $u_0=\arctan \sqrt{a-1}$ when $a>1$ are described by
\[\Im \{\psi(u)\}=\Im \{\psi(u_0)\}=-c,\]
where 
\[c=\sqrt{a-1}-a\,\arctan \sqrt{a-1} >0.\]
Setting $u=\xi+i\eta$, we have
\[\Im \frac{i(\tan \xi+i\tanh \eta)}{1-i\tan \xi \tanh \eta}-a\xi=\frac{\tan \xi(1-\tanh^2 \eta)}{1+\tan^2 \xi \tanh^2 \eta}-a\xi=-c,\]
which yields the steepest descent paths in the form
\bee\label{b1}
\tanh \eta=\pm \sqrt{\frac{\tan \xi-a\xi+c}{\tan \xi(1+(a\xi-c) \tan \xi)}}~.
\ee

Hence, $\eta=1$ when $(a\xi-c) \sec^2 \xi=0$; that is, when $\xi=c/a$. A similar argument applies to the saddle in the left-half plane to yield the result that the steepest descent paths $AB$ and $CD$ in Fig.~1(a) asymptote
to the line 
\bee\label{b2}
\Re (u)=\pm c/a.
\ee


\vspace{0.6cm}

\begin{thebibliography}{99}
\footnotesize{
\bibitem{AS}
M. Abramowitz and I.A. Stegun, {\it Handbook of Mathematical Functions}, U.S. National Bureau of Standards, Applied Mathematics Series, vol. 55, Washington, D.C., 1964.

\bibitem{ACM}
A. Apellblat, A. Consiglio and F. Mainardi, The Bateman functions revisited after 90 years: A survey of old and new results. MDPI Special Functions part III, Mathematics {\bf 9(11)} (2021) 1273. [arXiv:2104.08596].

\bibitem{HB}
H. Bateman, The $k$-function, a particular case of the confluent hypergeometric function, Trans. Amer. Math. Soc. {\bf 33} (1931) 817--831.

\bibitem{THH}
T.H. Havelock, The method of images in some problems of surface waves, Proc. Roy. Soc. London {\bf 108A} (1925) 582--591.

\bibitem{Olv}
F.W.J. Olver, {\it Asymptotics and Special Functions}, Academic Press, New York, 1974; Reprinted in A.K. Peters, Massachussets, 1997.

\bibitem{DLMF}
F.W.J. Olver, D.W. Lozier, R.F. Boisvert and C.W. Clark (eds.),    
{\it NIST Handbook of Mathematical Functions}, Cambridge University Press, Cambridge, 2010.

\bibitem{PHad}
R.B. Paris, {\it Hadamard Expansions and Hyperasymptotic Evaluation}, Encyclopedia of Mathematics and its Applications Vol. 141, Cambridge University Press, Cambridge, 2011.

}
\end{thebibliography}
\end{document}